\documentclass{article}
\usepackage[utf8]{inputenc}

\usepackage{amsmath}
\usepackage{amssymb}
\usepackage{type1cm}
\usepackage{amssymb}
\usepackage{amscd}
\usepackage{latexsym}
\usepackage{amsthm}

 \makeatletter
    
    \@addtoreset{equation}{section}
  \makeatother

\theoremstyle{definition}
\newtheorem{theo}{Theorem}[section]
\newtheorem{defi}[theo]{Definition}
\newtheorem{lem}[theo]{Lemma}
\newtheorem{prop}[theo]{Proposition}

\newtheorem{rem}[theo]{Remark}
\newtheorem{cor}[theo]{Corollary}
\newtheorem{ex}[theo]{Example}

\newtheorem{prf}{Proof}

\title{A small deformation and modifications of locally conformally balanced manifolds}
\author{Hirokazu Shimobe}
\date{}

\begin{document}

\maketitle

\begin{abstract}
We will consider locally conformally balanced manifolds. We prove that a locally conformally balanced condition is not stable under a small deformation. We prove that locally conformally balanced condition is stable under any proper modification. We prove that symmetric products of the Kodaira surface can be resolve to locally conformally balanced manifolds by Hilbert-Chow map. 
\end{abstract}

\section{Introduction}

Let $X$ be a complex manifold. A Hermitian metric on $X$ will be identified throughout this paper with the corresponding positive-definite $C^{\infty} (1,1)$-form $\omega$ on $X$. If $X$ admit a Hermitial metric $\omega$ which is closed, $X$ is called a K\"{a}hler manifold and  if $X$ admit a Hermitian metric $\omega$ which is co-closed, $X$ is called a balanced manifold. When one give a geometric structure, it is fundamental to consider if the given geometric structure is stable under proper modifications and small deformations. One can construct new compact complex manifolds with the same property if it dose. It is known that the K\"{a}hler condition is stable under any small deformation, but the balanced condition is not. For example, The Iwasawa manifold is balanced, but a small deformation of the Iwasawa manifold is not balanced. On the other hand, it is known that the K\"{a}hler condition is not stable under some proper modification \cite{Hironaka} and the balanced conditon is stable under any proper modification \cite{AB4}. 

We can consider metrics which locally conformal to given a metric, such as a K\"{a}hler metric and a balanced metric. We call complex manifolds with a metric locally conformal to a K\"{a}hler metric as locally conformally K\"{a}hler manifolds (LCK for short) and we define as follow for the balanced condition.

\begin{defi} \label{d1.1}
Let $(X, \omega)$ be a Hermitian manifold, $ n := $ dim $ X > 1 $, with $d \omega^{n-1} = \theta \wedge \omega^{n-1}$, where $\theta$ is a closed 1-form. Then $X$ is called a {\it locally conformally balanced} (LC-balanced for short) manifold. 
We call $\theta$ as a Lee form.
\end{defi}

As for LCK, a lot of research has been done so far. However, for LC-balanced, little research has been done. In this paper, we will consider locally conformally balanced manifolds.

As a first remark, we will note the following fact. That is, LCK manifolds are LC-balanced manifolds, but there are LC-balanced manifolds with no LCK structure. The most simple examples come form the following proposition. 

\begin{prop}
{\it Compact complex manifolds given by direct product of LC-balanced manifolds are LC-balanced.}
\end{prop}

\begin{prf}
See Proposition \ref{p3.1}.
\end{prf}

Therefore the direct product of LCK manifolds is LC-balanced. The direct product of LCK manifolds is not generally LCK. For example, the direct product of Vaisman manifolds is not LCK (see Corollary 3.3 of \cite{Ts}). Therefore the direct product of any submanifolds (dim $\ge 2$) of Hopf manifold give examples of LC-balanced manifolds that are not LCK.

Next, we will consider stability of the LC-balanced condition under small deformation. It is known that the LCK condition is not stable under a small deformation \cite{Belgun}. Using this result of \cite{Belgun}, we prove the following result.

\begin{prop} \label{p1.1}
{\it We follow the termes used in} \cite{Inoue} ({\it also see} \cite{Belgun}). {\it Let $X = S_{n;p,q,u}^{+}$ with $u \in \mathbb{C} \verb|\| \mathbb{R}$ be an Inoue surface of the type $\Gamma \verb|\| Sol_{1}^{'4}$. Then the direct product of $X$ cannot admit any LC-balanced metric.}
\end{prop}

\begin{prf}
See Proposition \ref{p3.2}.
\end{prf}

This proposition implies that the LC-balanced condition is not stable under a small deformation.

Next, we will consider whether the LC-balanced condition is stable under modification. It is know that the LCK condition is not stable under a proper modification \cite{Ornea}. For example, a blow up along submanifolds (dim $\ge$ 1) of a Vaisman manifold is not LCK. For LC-balanced manifolds, the following Theorem is known.

\begin{theo} \cite{S} \label{p1.2}
{\it Let $\mu \colon \tilde{X} \rightarrow X$ be a proper modification of compact complex manifolds. Then if $\tilde{X}$ is an LC-balanced manifold, $X$ is an LC-balanced manifold.}
\end{theo}

In this paper, we will prove the following theorem by using Theorem 2.3 of \cite{AB2}.

\begin{theo} \label{p1.3}
{\it Let X be a compact complex manifold endowed with a LC-balanced metric. Let $Y \subset X$ be a submanifold. Then the blow up $\tilde{X}$ of $X$ along $Y$ has a LC-balanced metric.} 
\end{theo}

\begin{prf}
See Theorem \ref{th4.1}
\end{prf}

By Theorem \ref{p1.2} and Theorem \ref{p1.3}, we can conclude that LC-balanced condition is stable under proper modification.

Finally, we will consider the Hilbert schemes of Kodaira surface. It is known that the Hilbert schemes of the Kodaira surface are non K\"{a}hler and are smooth holomorphic symplectice manifolds. We will study the geometric structure of these holomorphic symplectic manifolds. 

\begin{theo}
{\it Let $X$ be the Kodaira surface and $Y$ the symmetric product. Let $\tilde{Y}$ be a resolution of $Y$ by Hilbert-Chow map. Then $\tilde{Y}$ is LC-balanced.}
\end{theo}

\begin{prf}
See Theorem \ref{t5.1}
\end{prf}

That is, we obtain a lot of examples of holomorphic symplectic manifolds with LC-balanced structure.

\section{Preliminaries}

In this section, we give a brief introduction to locally conformally K\"{a}hler structures and locally conformally balanced structures. We define LCK manifolds as follow.

\begin{defi}
Let $(X, \omega)$ be a Hermitian manifold, $ n := $ dim $ X > 1 $, with $d \omega = \theta \wedge \omega$, where $\theta$ is a closed 1-form. Then $X$ is called a locally conformally K\"{a}hler (LCK for short) manifold. 
We call $\theta$ as a Lee form.
\end{defi}

A particular class of LCK manifolds are Vaisman manifolds.

\begin{defi}
Let $X$ be an LCK manifold and $\theta$ the Lee form. Then $X$ is called Vaisman manifold if $\nabla \theta = 0$ where $\nabla$ is the Levi-Civita connection.
\end{defi}

It is known that submanifolds of a Hopf manifold are Vaisman manifolds. 

It is trivial that submanifolds of LCK manifolds are LCK. It is known that an LCK condition is not stable under a small deformation and a proper modification \cite{Belgun} \cite{Ornea}. 

There is a conditon for a compact complex manifold to admit an LCK metric or an LC-balanced metric. This result was proved by Otiman \cite{O}.

\begin{prop} \cite{O} \label{p2.2}
{\it Let $X$ be a complex, compact manifold and $\theta$ a real closed 1-form.
There exists a transverse ($p$, $p$) $d_{\theta}$-closed form if and only if there are no positive currentswhich are ($p, p$)-components of $d_{\theta}$-boundaries.}
\end{prop}

LCK metrics and LC-balanced metrics can be defined on orbifolds. Let $X$ be an orbifold of dimension $n$. An orbifold chart on $X$ is a triple ($\tilde{U}, \Gamma, \phi$) where $\tilde{U}$ is a domain in $\mathbb{C}^{n}$, $\Gamma$ is a finite group acting effectively as automorphisms of $\tilde{U}$, and $\phi \colon \tilde{U} \rightarrow U$ is an analytic cover onto an open set $U \subset X$ such that $\phi \circ \gamma = \phi$ for every $\gamma \in \Gamma$ and induced natural map $\tilde{U} /\Gamma \rightarrow U$ is a homeomorphic. An orbifold atlas on $X$ is a family $\mathcal{U} = \{ (\tilde{U}_{i}, \Gamma_{i}, \phi_{i})\}$ of orbifold charts such that $X = \cup U_{i}$, where $U_{i} := \phi_{i} (\tilde{U}_{I})$.

We define a $(p,q)$-form on an orbifold $\mathcal{X} = (X, \mathcal{U})$ as follow. That is, a $(p,q)$-form on an orbifold $X$ is a collection $\tilde{\psi} = \{ \tilde{\psi}_{i} \}$ of smooth $(p,q)$-forms $\tilde{\psi}_{i}$ on $\tilde{U}_{i}$, for each orbifold chart $(\tilde{U}_{i}, \Gamma_{i}, \phi_{i})$, so that $\gamma^{*} \tilde{\psi_{i}} = \tilde{\psi_{i}}$ for each $\gamma \in \Gamma_{i}$, and $\lambda_{ji}^{*} \tilde{\psi_{j}} = \tilde{\psi_{i}}$ for each injection $\lambda_{ji} \colon \tilde{U}_{i} \rightarrow \tilde{U}_{j}$.

Therefore LCK structure and LC-balanced structure on an orbifold can be defined. 

A $(p,q)$-current on an orbifold $\mathcal{X} = (X, \mathcal{U})$ is defined as a collection $\tilde{T} = \{\tilde{T_{i}} \}$ of bidegree $(p,q)$-current $\tilde{T}_{i}$ on $\tilde{U}_{i}$, for each orbifold chart $(\tilde{U}_{i}, \Gamma_{i}, \phi_{i})$, so that $\gamma_{*} \tilde{T_{i}} = \tilde{T_{i}}$ for each $\gamma \in \Gamma_{i}$, and $(\lambda_{ji})_{*} \tilde{T_{i}} = \tilde{T_{j}}_{|\lambda_{ji}(\tilde{U_{i}})}$.

\begin{lem} \cite{C} \label{l2.1}
{\it Let $\mathcal{X} = (X, \mathcal{U})$ be an orbifold. Let $T$ be a positive closed bidegree (1,1) current. Then for any $x \in X$, there is a neighborhood $x \in U \subset X$ such that there is a plurisubharmonic (psh for short) function $v$ on $U$ with $d d^c v = T$ on $U \cap X_{reg}$.}
\end{lem}

\section{A small deformation of an LC-balanced manifold}

In this section, we will consider a small deformaiton of an LC-balanced manifold. First, we prove the following proposition.

\begin{prop} \label{p3.1}
{\it The direct product of LC-balanced manifolds are LC-balanced manifolds.}
\end{prop}

\begin{prf}
Let $X_{1}$ and $X_{2}$ be LC-balanced manifolds. Let $\omega_{i}$ be an LC-balanced metric and $\theta_{i}$ the Lee form of $X_{i}$, $i = 1,2$. Put that $n_{i} = $dim $X_{i}$, $i =1,2$. We put $\Omega := \omega_{1} + \omega_{2}$, and compute $d \Omega^{n_{1} + n_{2} - 1}$. Since 
\begin{equation*}
    \Omega^{n_{1} + n_{2} - 1} = c (\omega_{1}^{n_{1}} \wedge \omega_{2}^{n_{2} -1} + \omega_{1}^{n_{1} -1} \wedge \omega_{2}^{n_{2}} ),
\end{equation*}
where $c$ is a constant, 
\begin{align*}
    d \Omega^{n_{1} + n_{2} - 1} &= c (\omega_{1}^{n_{1}} \wedge \theta_{2} \wedge \omega_{2}^{n_{2} -1} + \theta_{1} \wedge \omega_{1}^{n_{1} -1} \wedge \omega_{2}^{n_{2}} ) \\
    &= (\theta_{1} + \theta_{2}) \wedge \Omega^{n_{1} + n_{2} - 1}.
\end{align*}
Therefore $\Omega$ is an LC-balanced metric on $X_{1} \times X_{2}$.\qed

\end{prf}

In general, the direct product of LCK manifolds is not LCK.

\begin{prop} \label{p3.1.1} \cite{Ts} 
{\it Let $X_{1}$ and $X_{2}$ be Vaisman manifolds. Then the direct product $X_{1} \times X_{2}$ is not LCK.}
\end{prop}

\begin{ex} \label{ex3.1}
Let $X$ and $Y$ be submanifolds of a Hopf manifold. Then by Proposition \ref{p3.1} and Proposition \ref{p3.1.1}, $X \times Y$ is an LC-balanced manifold that does not admit any LCK metrics.
\end{ex}

By Example \ref{ex3.1}, we get a lot of simple examples of LC-balanced manifolds. Using this example, we will discuss a small deformation of a LC-balanced manifold. Our proof of the following proposition is based on arguments of \cite{Belgun} (Proposition 18).

\begin{prop} \label{p3.2}
{\it {\it We follow the termes used in} \cite{Inoue} ({\it also see} \cite{Belgun}). {\it Let $X = S_{n;p,q,u}^{+}$ with $u \in \mathbb{C} \verb|\| \mathbb{R}$ be an Inoue surface of the type $\Gamma \verb|\| Sol_{1}^{'4}$. Then The direct product of $X$ cannot admit any LC-balanced metric.}}
\end{prop}

\begin{prf}
Let Lie($Sol_{1}^{'4}$) be a Lie algebra of $Sol_{1}^{'4}$ (see Proposition 17 of \cite{Belgun}). 

\begin{equation} \label{eq3.1}
    Lie(Sol_{1}^{'4}) = \langle Y, Z, T, U \rangle,
\end{equation}
where Z is central, $[Y, T] = Y, [T, U] = U, [Y, U] = Z$. 

Suppose that $X \times X$(or we index $X \times X$ as $X_1 \times X_2$) is LC-balanced. Let $\Omega$ be the LC-balanced metric and $\theta$ the Lee form.By Proposition 2.12 of \cite{Angella}, we may suppose that $\Omega$ and $\theta$ are left invariant. Since $\theta$ is $d$-closed, we can set that $\theta(T_i) =: k_i \in \mathbb{R}, \theta(Y_i) = \theta(Z_i) = \theta(U_i) = 0$, where $i = 1, 2$. Since $\Omega$ is a LC-balanced metric,
\begin{equation*}
    d \Omega^3 = \theta \wedge \Omega^3.
\end{equation*}
By definition and (\ref{eq3.1}), we compute the following derivation.
\begin{align*}
     d \Omega^3 (Y_1, Z_1, T_1, U_{1}, Y_2, Z_2, U_2) = c \{ \Omega^3(Y_1, Z_1, U_1, Y_2, Z_2, U_2) - \\ \Omega^3 (U_1, Y_1, Z_1, Y_2, Z_2, U_2) \}
     &= 0.
\end{align*}
On the other hand,
\begin{align*}
    \theta \wedge \Omega^3 (Y_1, Z_1, T_1, U_{1}, Y_2, Z_2, U_2) = \\
    k _{1} \Omega^3 (Y_1, Z_1, U_1, Y_2, Z_2, U_2).
\end{align*}
Therefore $k_1 = 0$. Similarly we have $k_2 = 0$. Therefore $\theta = 0$, but $X_1 \times X_2$ can not admit any balanced metric (if it were, the projection $X_1 \times X_2 \rightarrow X_1$ would imply that $X_1$ is balanced). Hence $X_1 \times X_2$ is not LC-balanced. \qed
\end{prf}


The direct product of the Inoue surface of Proposition \ref{p3.2} is a small deformation of an LC-balanced manifold. Therefore, the LC-balanced condition is not stable under a small deformation.

\begin{cor}
A small deformation of an LC-balanced manifold is not LC-balanced.
\end{cor}

\begin{rem}
It is known that a small deformation of a balanced manifold is not balanced. For example, the Iwasawa manifold is balanced, but a small deformation of the Iwasawa manifold is not balanced \cite{AB5}.
\end{rem}

\section{A proper modification of LC-balanced manifolds}
We will prove that any blow up of LC-balanced manifolds are LC-balanced. As a result , we conclude that LC-balanced condition is stable under any proper modification. 

\begin{theo} \label{th4.1}
{\it Let $X$ be an LC-balanced manifold and $\theta$ the Lee form. Let $f \colon \tilde{X} \rightarrow X$ be a blow up of $X$ along a submanifold $D$} ({\it codim $D \ge 2$}) {\it on $X$. Then $\tilde{X}$ is LC-balanced.} 
\end{theo}

\begin{prf}
Suppose that there exists a current $T = \partial_{\theta^{'}} \overline{S} + \overline{\partial}_{\theta^{'}} S$ on $\tilde{X}$, where $\theta^{'} = f^{*} \theta$. By Proposition \ref{p2.2}, we will prove that $T = 0$. Since $X$ is LC-balanced, Supp $ T \subset E := f^{-1} (D)$. In case that $D$ is a point, by \cite{L}, it is known that $\tilde{X}$ is LC-balanced. So we suppose that dim $D \ge 1$.  We divide $D$ into the following domain,

\begin{equation*}
    D = U_{1} \cup \cdots \cup U_{l},
\end{equation*}
where $j_{i}^{*} [\theta] = 0$ for embedding $j_{i} \colon U_{i} \rightarrow X$, $i = 1, \cdots, l$. Let $V_{i}$ be an open neighborhood of $U_{i}$ such that $H^{1} (V_{i}) \cong H^{1}(U_{i})$ and $V_{i} \cap D = U_{i}$. Since $\partial_{\theta^{'}} \overline{\partial_{\theta^{'}}} T = 0$, there exists a function $g_{i}$ on $f^{-1} (V_{i})$ such that $e^{g_{i}} \partial \overline{\partial} e^{-g_{i}} T = 0$ on $f^{-1} (V_{i}$). Let $T_{i}^{'} := e^{-g_{i}} T$. 

Since $T_{i}^{'}$is the component of a boundary,
\begin{equation} \label{eq4.1}
    T_{i}^{'} = \overline{\partial} S^{'} + \partial \overline{S^{'}}
\end{equation}
for a suitable current $S^{'}$ of bidegree (1,0). We use Theorem 2.3 of Alessandrini-Bassanelli \cite{AB2} that argued stability of balanced manifolds under modificaition.

\begin{theo} \cite{AB2} \label{th4.2}
{\it Let $M^{'}$ and $M$ be complex manifolds and $f \colon M^{'} \rightarrow M$ be the blow up of $M$ with smooth center $Y$. Let $T$ be a real $\partial \overline{\partial}$-closed current on $M^{'}$ of order zero and of bidegree} (1,1) {\it whose support is conttained in the exceptional set $E$.
Then there exists a pluriharmonic function $h \colon Y \rightarrow \mathbb{R}$ such that
\begin{equation*}
    T = (h \circ f) [E].
\end{equation*}
Moreover, if $T$ is the limit in the weak topology of currents that are components of boundaries, then T = 0.
\qed}
\end{theo}

$T^{'}$ in (\ref{eq4.1}) is a positive current, therefore, order zero. Hence since $T^{'}$ satisfy conditions of Theorem \ref{th4.2}, we can conclude that $T^{'} = 0$, that is, $T = 0$. 
\qed
\end{prf}

\begin{rem}
In Theorem \ref{th4.1}, if $D$ is a point or compact induced globally conformally balanced submanifold, it was proved that a blow up of $X$ along $D$ was LC-balanced \cite{L} \cite{Y}.
\end{rem}

It is known that any blow down of $X$ is LC-balanced~\cite{S}, so we have the following corollary.

\begin{cor}
{\it Let $\mu \colon \tilde{X} \rightarrow X$ be a proper modification of compact complex manifolds $X$ and $\tilde{X}$. Then $\tilde{X}$ is a LC-balanced manifold if and only if $X$ is a LC-balanced manifold.}
\end{cor}

\begin{ex}
Let $M$ be a Hopf manifold and $D$ a submanifold of $M$. We suppose that dim $D \ge 1$ and codim $D \ge 2$. Then a blow up $\tilde{M}$ of $M$ along $D$ is LC-balanced.
\end{ex}

\section{Holomorphic symplectic manifolds with LC-balanced structure}

Let X be an LCK surface and $X^{(n)}$ a symmetric product, that is, $X^{(n)} := X^n / S^n$, where $S^n$ is a symmetric group. Note that $X^{(n)}$ is an orbifold. There is Hilbert-Chow map $\pi \colon X^{[n]} \rightarrow X^{(n)}$ as a resolution of singlarities \cite{Fogarty}. Since dim $X = 2$, $X^{[n]}$ is a $2n$ dimensional compact complex manifold. When $X$ is the Kodaira surface, $X^{[n]}$ is a holomorphic symplectic manifold \cite{Beauville} \cite{Fujiki}. We will prove that $X^{[n]}$ admits an LC-balanced metric. 

Let $X_{*}^{(n)}$ be a subset of $X^{(n)}$ consisting of $\sum \nu_{i}[x_i]$ ($x_i$:distinct) with $\nu \leq 2, \nu_2 = \cdots \nu_k = 1$. Let $pr \colon X^n \rightarrow X^{(n)}$ be a projection. Define $X_{*}^{[n]}$ and $X_{*}^{n}$ as follows.

\begin{align*}
    X_{*}^{[n]} := \pi^{-1} (X_{*}^{(n)}) \\
    X_{*}^{n} := pr^{-1}(X_{*}^{(n)}).
\end{align*}
Let us denote by $\Delta \subset X^n$ the “big diagonal” consisting of elements ($x_1,\dots,x_n$) with $x_i = x_j$ for some $i = j$. Then we get the following commutative diagram.

$$
\begin{CD}
Blow_{\Delta} (X_{*}^{n}) @>\eta >> X_{*}^{n} \\
@V \rho VV @VVV \\
X_{*}^{[n]} @>\pi>> X_{*}^{(n)},
\end{CD}
$$
where $\eta$ is the blow up along $\Delta$ of $X_{*}^{n}$ and $\rho$ is the map given by taking the quotient by the action of $S^{n}$. 
\begin{theo} \label{t5.1}
{\it Let $X$ be the Kodaira surface. Let $X^{[n]}$ be a holomorphic symplectic manifold as above. Then $X^{[n]}$ admit an LC-balanced metric.}
\end{theo}

\begin{prf}
Let $\theta$ be a Lee form of $X$. Since $X$ is a nilmanifold, we may suppose that $\theta$ is an invariant closed form. By Proposition \ref{p3.1}, $X^n$ is an LC-balanced manifold. Let $\theta_{1} + \cdots + \theta_{n}$ be the Lee form, where $\theta_{i} = \theta$ on $X$. Since $\pi$ is Hilbert-Chow resolution, the pull-back $\theta^{'} = \pi^{*} (\theta_{1} + \cdots + \theta_{n})$ can be defined. We may assume that $\theta^{'}$ can be assumed to be smooth on $X^{[n]}$ (One can be proven by the same method that a holomorphic symplectic form can be defined on $X^{[n]}$).

Suppose that $X^{[n]}$ is not LC-balanced. Then there is a bidegree (1,1) positive current $T = \partial_{\theta^{'}} \overline{S} + \overline{\partial_{\theta^{'}}} S$ \cite{O}. Since $X^{n}$ is LC-balanced, Supp $T \subset E$, where $E$ is exceptional divisor of $\pi$. Let $E = A \cup B$, where $A$ is the union of irreducible components of $E$ of codimension bigger than 1, $B = \cup_{i} E_{i}$ is the union of irreducible components of $E$ of codimention 1.

Let $D_{i} := E_{i} \cap X_{*}^{[n]}$, then $D_{i}$ is smooth. We divide $D_{i}$ into the following domain,

\begin{equation*}
    D_{i} = U_{1} \cup \cdots \cup U_{l},
\end{equation*}
where $j_{k}^{*} [\theta] = 0$ for embedding $j_{k} \colon U_{k} \rightarrow X^{[n]}$, $k = 1, \cdots, l$. Let $V_{k}$ be an open neighborhood of $U_{k}$ such that $H^{1} (V_{k}) \cong H^{1}(U_{k})$ and $V_{k} \cap D_{i} = U_{k}$. We have 
\begin{equation*}
    e^{g_{k}} \partial \overline{\partial} e^{-g_{k}} T = 0
\end{equation*}
on $V_{k}$. Let $T_{k}^{'} := e^{-g_{k}} T$ on $V_{k}$. $T_{k}^{'}$ of (\ref{eq4.1}) is a positive current, therefore, order zero. Hence since $T_{k}^{'}$ satisfy conditions of Theorem \ref{th4.2}, we can conclude that $T_{k}^{'} = 0$, that is, Supp$T \subset A$. Since codim $A > 1$, we can conclude $T = 0$. Therefore $X^{[n]}$ admit a LC-balanced metric. 

\end{prf}

Finally, we prove that $X^{[n]}$ does not admit any balanced metrics. 
\begin{prop}
{\it Let $X^{[n]}$ be a holomorphic symplectic manifold as above. Then $X^{[n]}$ is not balanced.}
\end{prop}

\begin{prf}
We follow termes and arguments of \cite{P2}. We will prove that $X^{[n]}$ is not strongly Gauduchon (sG for short), therefore $X^{[n]}$ is not balanced. Suppose that $X^{[n]}$ is sG. Then we should prove that $X^{n}$ is sG. Suppose that $X^{n}$ is not sG. Then by Proposition 3.3 in \cite{P1}, there is no non-zero current $T$ of bidegree (1, 1) on $X^n$ such that $T \ge 0$ and $T$ is $d$-exact on $X^n$. Let $p \colon X^n \rightarrow X^{(n)}$ be a projection. Then $\tilde{T} := p_{*} T$ is well-defined current on $X^{(n)}$ (Note that $X^{(n)}$ is a orbifold). $\tilde{T}$ is a $d$-exact non zero positive (1,1) current. By Lemma \ref{l2.1}, for any $x \in X^{(n)}$, there is a neighborhood $U \subset X^{(n)}$ of $x$ such that $T|_{U \cap X_{reg}^{(n)}} = i \partial \overline{\partial} v$, where $v$ is a psh on $U$. Define that $(\pi^{*} \tilde{T})|_{\pi^{-1} (U)} := i \partial \overline{\partial} (v \circ \pi)$. Then we can construct $\pi^{*} \tilde{T}$ as pull-back of $T$ by $\pi$. We will prove that $\pi^{*} T$ is $d$-exact on $X^{[n]}$.

By Main Theorem 1.1 of \cite{Dem}, there is a $\tilde{v_{j}} \in [\tilde{T}]_{DR}$ such that $\tilde{v_{j}} \xrightarrow{w} \tilde{T}$ and $\tilde{v_{j}} \ge -c \omega$, where $c$ is a positive constant and $j \in \mathbb{N}$. Since $\tilde{T}$ is $d$-exact, $\tilde{v_{j}}$ is $d$-exact. 

We will prove that $\pi^{*} \tilde{v_j} = d (\pi^{*} \tilde{u_j}) \xrightarrow{w} \pi^{*} \tilde{T}$. The argument is virtually the same as that of \cite{P2}. Hence we give a sketch of proof.

Pick smooth invariant form $\alpha^{\Gamma}$ ($X^{(n)}$ is an orbifold). 

\begin{equation*}
    \tilde{v_{j}} = \alpha^{\Gamma} + i \partial \overline{\partial} \psi_{j} \ge - C \omega,
\end{equation*}
where we choose a smooth function on $X^{2n}$, $\psi_{j}$, such that $\int_{X^n} \psi_{j} \omega^{2n} = 0$, where $\omega$ is a Hermitian metric on $X^{2n}$. By the defintion of Laplacian $\Delta$, we have
\begin{equation*}
    \Delta_{\omega} \psi_{j} = {\rm Trace}_{\omega} (\tilde{v_{j}} - \alpha).
\end{equation*}
Let $G$ be Grenn operator of $\Delta_{\omega}$,
\begin{equation*}
    \psi_{j} = G {\rm Trace}_{\omega} (\tilde{v_j} - \alpha).
\end{equation*}
Since $G$ is a compact operator from the Banach space of bounded Borel measures on $X^n$ to $L^1(X^n)$ and since the forms $\tilde{v_j}$ converge weakly to $T$, we
infer that some subsequence ($\psi_{j_k}$)$_k$ converges to a limit $\psi \in L^1(X^n)^\Gamma$ in $L^1(X^n)$-topology. Therefore, by the weak continuity of $\partial \overline{\partial}$, we have
\begin{equation*}
    \tilde{T} = \lim_{k} (\alpha + i \partial \overline{\partial} \psi_{j_{k}})
    = \alpha + i \partial \overline{\partial} \psi.
\end{equation*}
Since $(\psi_{j})_{j}$ is uniformly bounded above, $(\psi_{j} \circ \pi)_{j}$ is also uniformly bounded above on $X^{[n]}$. On the other hand, $\psi_{j_{k}} \circ \pi$ converges almost everywhere to $\psi \circ \pi$ on $X^{[n]}$. Therefore,
\begin{equation*}
    \pi^{*} \tilde{v_{j_{k}}} = \pi^{*} \alpha + i \partial \overline{\partial}(\psi_{j_{k}} \circ \pi)
    \xrightarrow{w} \pi^{*} \tilde{T} = \pi^{*} \alpha + \partial \overline{\partial} (\psi \circ \pi).
\end{equation*}
since $\pi^{*} \tilde{T}$ is a non-zero positive $(1,1)$ $d$-exact current, $X^{[n]}$ is non-sG. That is, $X^{[n]}$ is not balanced.
\qed
\end{prf}

E-mail:rsc96831@gmail.com
\end{document}